\documentclass[10pt,leqno]{amsart} 

   \usepackage[centertags]{amsmath}
   \usepackage{amsfonts}
   \usepackage{amsmath}
   \usepackage{amssymb}
   \usepackage{amsthm}
   \usepackage{newlfont}
   \usepackage[latin1]{inputenc}
\usepackage{xcolor}

\allowdisplaybreaks[3]

\theoremstyle{plain}
\newtheorem{Theo}{Theorem}[section]

\newtheorem{Cor}[Theo]{Corollary}
\newtheorem{Lem}[Theo]{Lemma}

\theoremstyle{definition}
\newtheorem{Def}{Definition}[section]

\theoremstyle{remark}
\newtheorem{Rem}{Remark}

\newcommand{\RR}{\mathbb{R}}
\newcommand{\ZZ}{\mathbb{Z}}

\newcommand{\co}{\operatorname{co}}
\newcommand{\He}{\operatorname{He}}

\numberwithin{equation}{section}

\newcommand{\card}{\operatorname{card}}

\newcommand\supp{\operatorname{supp}}

\newcommand\sign{\operatorname{sign}}

\newcommand*\pFq{\begingroup
        \dopFq
}
\catcode`\,12
\def\dopFq#1#2#3#4#5{%
        {}_{#1}F_{#2}\biggl(\genfrac..{0pt}{}{#3}{#4};#5\biggr)%
        \endgroup
}

\title[Zeros of linear combinations of orthogonal polynomials]{Zeros of linear combinations of orthogonal polynomials}
\author{Antonio J. Dur\'an}
\address{Departamento de An\'a\-li\-sis Mate\-m\'a\-ti\-co and IMUS,
        Universidad de Sevilla,
        41080 Sevilla, Spain}
\email{duran@us.es}
   \date{}

   \thanks{This research was partially supported by PID2021-124332NB-C21
(Minis\-te\-rio de Cien\-cia e Inno\-va\-ci\'on and Feder Funds (European Union)), and
FQM-262 (Jun\-ta de Anda\-lu\-c\'ia).}



\begin{document}
   \maketitle

\begin{abstract}
Given a sequence of orthogonal polynomials $(p_n)_n$ with respect to a positive measure in the real line, we study the real zeros of finite combinations of $K+1$ consecutive orthogonal polynomials of the form
$$
q_n(x)=\sum_{j=0}^K\gamma_jp_{n-j}(x),\quad n\ge K,
$$
where $\gamma_j$, $j=0,\cdots ,K$, are real numbers with $\gamma_0=1$, $\gamma_K\not =0$ (which do not depend on $n$).
We prove that for every positive measure $\mu$ there always exists a sequence of orthogonal polynomials with respect to $\mu$ such that
all the zeros of the polynomial $q_n$ above are real and simple for $n\ge n_0$, where $n_0$ is a positive integer
depending on $K$ and the $\gamma_j$'s.

\end{abstract}

\section{Introduction}
Given a sequence of orthogonal polynomials $(p_n)_n$ with respect to a positive measure in the real line (and infinitely many points in its support), we fix a positive integer $K$ and $K+1$ real numbers
$\gamma_j$, $j=0,\cdots ,K$, with $\gamma_0=1$, $\gamma_K\not =0$, and consider the sequence of polynomials
\begin{equation}\label{dqnv}
q_n(x)=\sum_{j=0}^K\gamma_jp_{n-j}(x),\quad n\ge K.
\end{equation}
Shohat \cite{Sho} was probably the first to observe that the orthogonality of the sequence $(p_n)_n$ implies that $q_n$ has at least $n-K$ real zeros in the convex hull of the support of $\mu$, $\co (\supp \mu)$, using the usual proof that $p_n$ has its $n$ zeros in $\co (\supp \mu)$ (see Lemma \ref{qor} below). Fejer \cite{Fej} had already studied when all the zeros of $q_n$ are real for the particular case $K=2$ in connection with quadrature formulas for the measure $\mu$, and Erdös and Turán \cite{ErTu} had also got some related results. More recent results can be found in \cite{Peh1,Peh2,BDr,DJM,KRZ} or \cite{IsSa,IsNo,INS} (for the case $n=K$).

The purpose of this paper is to go beyond the well-known and classic Shohat result in the study of the real zeros of the polynomials $q_n$, $n\ge K$.
We have to notice that this problem is strongly dependent on the normalization of the orthogonal polynomials $(p_n)_n$. Along this paper we denote by $(\hat p_n)_n$ the sequence of monic orthogonal polynomials.

As the main result, we prove in Section  \ref{sec1} the following unexpected result: we can always renormalize
the monic polynomials with a normalization sequence $(\rho_n)_n$, $p_n(x)=\rho_n\hat p_n(x)$, in such a way that given any positive integer $K$ and $K+1$ real numbers
$\gamma_j$, $j=0,\cdots ,K$, with $\gamma_0=1$, $\gamma_K\not =0$, there exists $n_0$ (depending on $K$ and the $\gamma_j$'s) such that all the zeros of the polynomial $q_n$ (\ref{dqnv}) are real and simple for $n\ge n_0$
and interlace the zeros of $p_{n-1}$ (for the definition of the interlacing property see Definition \ref{dip} below). We also provide explicit estimates for the normalization sequence $(\rho_n)_n$ (and then for the positive integer $n_0$) to satisfy that property (see Theorem \ref{desco} and Corollary \ref{cop}, bellow). We also discuss necessary and sufficient conditions for the interlacing properties of the zeros of $q_{n+1}$ and $q_n$.

We finally study some significant examples as the classical orthogonal polynomials, showing that our result applies to the Hermite and probabilistic Hermite polynomials, monic Laguerre polynomials and the normalization of the Jacobi polynomials given by $n!P^{\alpha,\beta}_n$.

\section{Preliminaries}
Let $(\hat p_n)_n$ be the monic sequence of orthogonal polynomials with respect to a positive measure $\mu$ supported in the real line (and with infinitely many points in its support).
They then satisfy a three term recurrence relation of the form
\begin{equation}\label{ttrrm}
\mbox{$x\hat p_n(x)=\hat p_{n+1}(x)+b_n\hat p_n(x)+\hat c_n\hat p_{n-1}(x)$, $n\ge 1$, with $\hat c_n>0, n\ge 1$.}
\end{equation}
If we re-normalize the monic family using a sequence of positive real numbers $(\rho_n)_n$, $\rho_0=1$, defining
$$
p_n(x)=\rho_n\hat p_n(x),
$$
then the sequence $(p_n)_n$ satisfies
the three term recurrence relation
\begin{equation}\label{ttrr}
xp_n(x)=a_{n+1}p_{n+1}(x)+b_np_n(x)+c_np_{n-1}(x),\quad n\ge 1,
\end{equation}
where
\begin{equation}\label{rmn}
a_{n}=\frac{\rho_{n-1}}{\rho_n},\quad c_n=\frac{\rho_{n}}{\rho_{n-1}}\hat c_n>0,\quad n\ge 1.
\end{equation}

Along this paper, the interlacing property  is defined as follows.

\begin{Def}\label{dip} Given two finite sets $U$ and $V$ of real numbers ordered by size, we say that $U$ interlaces $V$ if $\min U<\min V$ and between any two consecutive elements of any of the two sets there exists one element of the other.
\end{Def}
Observe that if $U$ interlaces $V$, then either $\card(U) = \card(V)$, and then $\max U < \max V$, or $\card(U) = 1 + \card(V)$, and then $\max U > \max V$. Observe also that the interlacing property is not symmetric, due to the condition $\min U < \min V$.

We will use the following version of Obreshkov Theorem (see \cite[Th. 8]{Branden}).

\begin{Theo}\label{obre}
Let $p$ and $q$ be real polynomials with $\deg p=1+\deg q$ and assume that $p$ and $q$ have real and simple zeros.
Then the zeros of $p$ interlace the zeros of $q$ if and only if for all real $\lambda $ the polynomial $p(z)+\lambda  q(z)$ has only real and simple zeros.
\end{Theo}

Along this paper we use the following notation: given $K+1$ real numbers $\gamma_0,\cdots ,\gamma_K$, we denote
\begin{equation}\label{dgam}
\Gamma_i=\max\{|\gamma_j|: i\le j\le K\}.
\end{equation}

The following elementary lemma will be useful.

\begin{Lem}\label{lem2} Define from the numbers $A_j$, $j=0,\cdots , K$, $A_0,A_K\not =0$, the polynomial $P_A$ as
$$
P_A(x)=\sum_{j=0}^KA_j x^{K-j}.
$$
If  $\theta $ is a zero of $P_A$, we define the polynomial $P_{B}$ and the numbers $B_j$, $j=0,\cdots , K-1$, as
$$
P_{B}(x)=\frac{P_{A}(x)}{x-\theta}=\sum_{j=0}^{K-1}B_j x^{K-1-j}.
$$
Then, on the one hand, we have
\begin{equation}\label{id1l2}
A_j=\begin{cases} B_j-\theta B_{j-1},& j=1,\cdots, K-1,\\ B_0,& j=0,\\
-\theta B_{K-1},&j=K.\end{cases}
\end{equation}
And on the other hand
\begin{equation}\label{id2l2}
B_{j}=\sum_{i=0}^{j}\theta^iA_{j-i},\quad 0\le j\le K-1.
\end{equation}
\end{Lem}

\section{Zeros of linear combinations of orthogonal polynomials}\label{sec1}
If $(p_n)_n$ is a sequence of orthogonal polynomials with respect to a positive measure $\mu$ supported in the real line (and with infinitely many points in its support),
Shohat \cite{Sho} proved that for a positive integer $K$, any finite combination of the consecutive orthogonal polynomials $p_n,\cdots ,p_{n-K}$ has always at least $n-K$ real zeros. Moreover, they are in the convex hull of the support of the measure.

For the sake of completeness we include a proof of this fact.

\begin{Lem}\label{qor} For real numbers $\gamma_j$, $j=0,\cdots, K$, $\gamma_0=1$, $\gamma_K\not =0$, define
$q(x)=\sum_{j=0}^K \gamma_jp_{n-j}(x)$, $n\ge K$. Then $q$ has at least $n-K$ real zeros in the convex hull of the support of the measure. Moreover, if $q$ has exactly $n-K$ real zeros, then they are simple.
\end{Lem}

\begin{proof}
Let $\xi_i$ be a real zero of $q$ of odd multiplicity belonging to the convex hull $S$ of the support of the measure $\mu$. Consider then the polynomial $r(x)=\prod_{i=1}^N(x-\xi_i)$. On the one hand, it is clear that $r(x)q(x)$ has constant sign in $S$. On the other hand, if $N<n-K$ then, because of the orthogonality of the polynomials $(p_n)_n$, we would have
$$
\int _S r(x)q(x)d\mu=0.
$$
Hence, we deduce that $N\ge n-K$. This proves the lemma.
\end{proof}

In this section, we improve the previous lemma. Indeed, we prove that there are normalization sequences $(\rho_n)_n$ such that the normalized orthogonal polynomials $p_n(x)=\rho_n\hat p_n(x)$ enjoy the following property: for any positive integer $K$, any finite combination
\begin{equation}\label{dqn}
q_n(x)=\sum_{j=0}^K\gamma_jp_{n-j}(x),\quad n\ge K,\quad \gamma_0,\gamma_K\not =0, \gamma_j\in \RR,
\end{equation}
has only real zeros for $n$ big enough (depending on the coefficients $\gamma_j$'s). Our proof is constructive, in the sense that it will
allow us to estimate the size of the normalization sequence $(\rho_n)_n$ and how big $n$ has to be to guarantee the real-rootedness of the polynomials $q_n$.

The proof is based in the following lemma.

\begin{Lem} Let $(p_n)_n$ be a sequence of orthogonal polynomials with respect to a positive measure satisfying (\ref{ttrr}).
Then, for real numbers $\gamma_j$, $j=0,\cdots, K$, with $\gamma_0,\gamma_K\not=0$,
there exist polynomials $A$, of degree $K-2$, and $B$, of degree $K-1$, such that the polynomials $q_n$ (\ref{dqn}) can be written as
\begin{equation}\label{dec}
q_n(x)=A(x)p_n(x)+B(x)p_{n-1}(x).
\end{equation}
Moreover, assume that $(\tau_n)_n$ is a decreasing sequence of positive numbers such that for certain positive integer $n_1$
\begin{equation}\label{cttrr}
\max\left\{\left|\frac{a_n}{c_{n-1}}\right|,\left|\frac{b_{n-1}}{c_{n-1}}\right|,\left|\frac{1}{c_{n-1}}\right|\right\}\le \tau_n\le 1/2,\quad n\ge n_1.
\end{equation}
Then for $n\ge n_1$, we have
\begin{align}\label{cotA}
|A(x)|&\le 2^K\Gamma_0\left(1+\sum_{j=1}^{K-2}\tau_{n-K+1}^{j+1}|x|^j\right),\\\label{cotB}
|B(x)|&\le 2^K\Gamma_1\sum_{j=0}^{K-1}\tau_{n-K+1}^{j}|x|^j,
\end{align}
where $\Gamma_i$ is defined by (\ref{dgam}).
\end{Lem}

\begin{proof}
Let us remark that we  do not assume here the normalization $\gamma_0=1$ for the polynomials $q_n$.

We proceed by induction on $K\ge 2$.

The case $K=2$ is straightforward. Indeed, using the recurrence relation (\ref{ttrr}), we can rewrite
$$
q_n(x)=\left(\gamma_0-\gamma_2\frac{a_n}{c_{n-1}}\right)p_n(x)
+\left(\gamma_1+\gamma_2\frac{x-b_{n-1}}{c_{n-1}}\right)p_{n-1}(x).
$$
So that,
$$
A(x)=\gamma_0-\gamma_2\frac{a_n}{c_{n-1}},\quad B(x)=\gamma_1+\gamma_2\frac{x-b_{n-1}}{c_{n-1}},
$$
from where the estimates (\ref{cotA}) and (\ref{cotB}) follow easily using (\ref{cttrr}).

Assume next that the lemma is true for $K-1$.
If we set
$$
\tilde q_{n-1}(x)=\sum_{j=1}^K\gamma_jp_{n-j}(x)=\sum_{j=0}^{K-1}\gamma_{j+1}p_{n-1-j}(x),
$$
we can write $
q_n(x)=\gamma_0 p_n(x)+\tilde q_{n-1}(x)$.

Using the induction hypothesis, we have
$$
\tilde q_{n-1}(x)=\tilde A(x)p_{n-1}(x)+\tilde B(x)p_{n-2}(x),
$$
where the polynomials $\tilde A$ and $\tilde B$ are associated to the real numbers $\tilde \gamma_{i}=\gamma_{i+1}$, $i=0,\dots, K-1$, and $K\to K-1$, $n\to n-1$. Hence, using (\ref{cotA}) and (\ref{cotB}), we deduce that $\tilde A$ and $\tilde B$ satisfy
\begin{align}\label{cotAx}
|\tilde A(x)|&\le 2^{K-1}\Gamma_1\left(1+\sum_{j=1}^{K-3}\tau_{n-K+1}^{j+1}|x|^j\right),\\\label{cotBx}
|\tilde B(x)|&\le 2^{K-1}\Gamma_2\sum_{j=0}^{K-2}\tau_{n-K+1}^{j}|x|^j.
\end{align}

Using again the recurrence relation (\ref{ttrr}), we can rewrite
$$
q_n(x)=\left(\gamma_0-\frac{a_n}{c_{n-1}}\tilde B (x)\right)p_n(x)
+\left(\tilde A(x)+\frac{x-b_{n-1}}{c_{n-1}}\tilde B(x)\right)p_{n-1}(x).
$$
This gives
\begin{equation}\label{decg}
A(x)=\gamma_0-\frac{a_n}{c_{n-1}}\tilde B(x),\quad B(x)=\tilde A(x)+\frac{x-b_{n-1}}{c_{n-1}}\tilde B(x).
\end{equation}
On the one hand,  using (\ref{cttrr}) and the bounds (\ref{cotBx}) for $\tilde B$, we deduce
$$
|A(x)|\le |\gamma_0|+\tau_n 2^{K-1}\Gamma_2\sum_{j=0}^{K-2}\tau_{n-K+1}^{j}|x|^j\le 2^K\Gamma_0\left(1+\sum_{j=1}^{K-2}\tau_{n-K+1}^{j+1}|x|^j\right)
$$
On the other hand, using (\ref{cttrr}) and the bounds (\ref{cotAx}) for $\tilde A$ and (\ref{cotBx}) for $\tilde B$, we deduce
\begin{align*}
|B(x)|&\le |\tilde A(x)|+\tau_n|\tilde B(x)|+\tau_n |x\tilde B(x)|\\
&\le 2^{K-1}\Gamma_1\left[1+\sum_{j=1}^{K-3}\tau_{n-K+1}^{j+1}|x|^j+\sum_{j=0}^{K-2}\tau_{n-K+1}^{j+1}|x|^j
+\sum_{j=0}^{K-2}\tau_{n-K+1}^{j+1}|x|^{j+1}\right]\\
&\le 2^K\Gamma_1\left(\sum_{j=0}^{K-1}\tau_{n-K+1}^{j}|x|^j\right).
\end{align*}
\end{proof}

We are now ready to prove the main result of this section.

\begin{Theo}\label{desco} Let $(\hat p_n)_n$ be a sequence of monic orthogonal polynomials with respect to a positive measure satisfying (\ref{ttrrm}). Write $\zeta_{j,n}$, $j=1,\dots, n$, for the $n$ real zeros of $\hat p_n$
arranged in increasing order. Given a sequence of numbers $(\rho_n)_n$, $\rho_n\not=0$, $n\ge0$, assume that there exist a positive constant $c$, a positive integer $n_1$ and a decreasing sequence of positive numbers $(\tau_n)_n$ such that $\lim_n\tau_n=0$ and
\begin{align}\label{c3t2}
\max\left\{\left|\frac{\rho_{n-2}}{\rho_{n}\hat c_{n-1}}\right|,\left|\frac{\rho_{n-2}b_{n-1}}{\rho_{n-1}\hat c_{n-1}}\right|,\left|\frac{\rho_{n-2}}{\rho_{n-1}\hat c_{n-1}}\right|\right\}&\le \tau_n,\quad n\ge n_1\\\label{c3t3}
\max\{|\zeta_{j,s}|:1\le j\le s, n-1\le s\le 2n\}&\le \frac{c}{\tau_n},\quad n\ge n_1.
\end{align}
Define next the sequence of normalized orthogonal polynomials $p_n(x)=\rho_n\hat p_n(x)$.
For any positive integer $K$ and any finite set of $K+1$ real numbers $\gamma_j$, $j=0,\dots, K$, with $\gamma_0=1$ and $\gamma_K\not=0$, let $n_0$ be the first positive integer such that
\begin{equation}\label{dta}
\mbox{$n_0\ge \max\{2K,n_1\}$ \quad and  \quad $\displaystyle \tau_{n_0}<\min \left\{\frac12, \frac{1}{2^{K-1}\Gamma_2\sum_{j=0}^{K-2}c^j}\right\} $.}
\end{equation}
Then, for $n\ge n_0$ the polynomial $q_n$ (\ref{dqn}) has only real zeros. Moreover the zeros are simple and interlace the zeros of $p_{n-1}$.
\end{Theo}

\begin{proof}
Since the polynomials $(p_n)_n$ satisfy the recurrence relation (\ref{ttrr}), the identities (\ref{rmn}) and (\ref{c3t2}) show that the sequences $(a_n)_n$, $(b_n)_n$, $(c_n)_n$ and $(\tau_n)_n$ satisfy (\ref{cttrr}) for $n\ge n_0$. Hence, there exist polynomials $A$ and $B$ such that
\begin{equation}\label{expi}
q_n(x)=A(x)p_{n}(x)+B(x)p_{n-1}(x),\quad n\ge 1,
\end{equation}
and $A$ and $B$ satisfy the bounds (\ref{cotA}) and (\ref{cotB}).

We next prove that, for $n\ge n_0$,
$$
A(\zeta_{l,n-1})>0,\quad l=1,\dots, n-1,
$$
where $\zeta_{l,n-1}$, $l=1,\dots, n-1$, are the $n-1$ real zeros of $\hat p_{n-1}$ arranged in increasing order.
Indeed, using the first identity in (\ref{decg}), (\ref{cotBx}) for $\tilde B$, and (\ref{cttrr}), we have
\begin{align*}
A(\zeta_{l,n-1})&\ge 1-\left|\frac{a_n}{c_{n-1}}\tilde B(\zeta_{l,n-1})\right|\\ & \ge 1-\tau_n 2^{K-1}\Gamma_2\sum_{j=0}^{K-2}\tau_{n-K+1}^{j}|\zeta_{l,n-1}|^j.
\end{align*}
Using (\ref{c3t3}), we have for $n\ge n_0\ge 2K$ that $n-K\le n-1\le 2n-2K+2$, and hence $\tau_{n-K+1}|\zeta_{l,n-1}|\le c$. We then have from (\ref{dta})
\begin{align*}
A(\zeta_{l,n-1})&\ge 1-\tau_n 2^{K-1}\Gamma_2\sum_{j=0}^{K-2}\tau_{n-K+1}^{j}|\zeta_{l,n-1}|^j\\
&\ge 1-\tau_n 2^{K-1}\Gamma_2\sum_{j=0}^{K-2}c^j \\
&\ge 1-\tau_{n_0}2^{K-1}\Gamma_2\sum_{j=0}^{K-2}c^j \\
&> 0.
\end{align*}

Hence, (\ref{expi}) gives
$$
q_n(\zeta_{l,n-1})q_n(\zeta_{l+1,n-1})=A(\zeta_{l,n-1})A(\zeta_{l+1,n-1})p_n(\zeta_{l,n-1})p_n(\zeta_{l+1,n-1}),
$$
and so
$$
\sign(q_n(\zeta_{l,n-1})q_n(\zeta_{l+1,n-1}))=\sign(p_n(\zeta_{l,n-1})p_n(\zeta_{l+1,n-1}))<0,
$$
because $p_n$ and $p_{n-1}$ are orthogonal polynomials and then the zeros of $p_{n}$ interlace the zeros of $p_{n-1}$.

This implies that $q_n$ has at least one zero in each interval $(\zeta_{l,n-1},\zeta_{l+1,n-1})$, $l=1,\dots,n-1$. Since the leading coefficients of $q_n$ and $p_n$ are the same, we deduce that $q_n$ has also at least one zero in the intervals $(-\infty,\zeta_{1,n-1})$ and $(\zeta_{n-1,n-1},+\infty)$. That is, $q_n$ has $n$ real and simple zeros which interlace the zeros of $p_{n-1}$.
\end{proof}

Theorem \ref{desco} has the following surprising consequence.

\begin{Cor}\label{cop} Let $(\hat p_n)_n$ be a sequence of monic orthogonal polynomials with respect to a positive measure satisfying (\ref{ttrrm}). Then, there exist sequences of positive numbers $(\rho_n)_n$ such that the sequence of normalized orthogonal polynomials $p_n(x)=\rho_n\hat p_n(x)$ satisfies the following: for any positive integer $K$ and any finite set of $K+1$ real numbers $\gamma_j$, $j=0,\dots, K$, with $\gamma_0=1$ and $\gamma_K\not=0$, the polynomials $q_n$ (\ref{dqn}) have only real zeros for $n$ big enough (depending only on $K$ and $\gamma_j$, $j=0,\dots, K$).
\end{Cor}

\begin{proof}
Write as before $\zeta_{j,n}$, $j=1,\dots, n$, for the $n$ real zeros of $\hat p_n$
arranged in increasing order.

We fix a decreasing sequence $(\theta_n)_n$ of positive numbers such that $\theta_n\le 1/2$ and $\lim_n\theta_n=0$. We next define the sequences $(\tau_n)_n$, $(e_n)_n$ and $(d_n)_n$ of positive numbers as follows: $\tau_0=1/2$,
\begin{align*}
\tau_n&=\min\{\theta_n,\tau_{n-1},1/\varsigma_{s} : n-1\le s\le 2n\},\quad n\ge 1,
\end{align*}
where $\varsigma_{s}=\max\{|\zeta_{j,s}|:1\le j\le s\}>0$, and
\begin{align*}
e_n&=\begin{cases} \frac{|\hat c_{n-1}|}{|b_{n-1}|}\tau_n,&b_{n-1}\not=0,\\ 1,&b_{n-1}=0,\end{cases}\\
d_n&=\min \{ |\hat c_{n-1}|\tau_n,e_n,9/10\}.
\end{align*}

We finally define the sequence $(\rho_n)_n$ of positive numbers by $\rho_0=1$ and
$$
\rho_n=\frac{\rho_{n-1}}{d_{n+1}}, \quad n\ge 1.
$$
It is enough then to prove that the sequences $(\rho_n)_n$ and $(\tau_n)_n$ satisfy the hypothesis of Theorem \ref{desco}.

Indeed, since $0<d_n< 1$, we deduce that $(\rho_n)_n$ is an increasing sequence of positive numbers. By definition, $(\tau_n)_n$ is a decreasing sequence of positive numbers with $\tau_n\le 1/2$ and $\lim_n\tau_n=0$. Since
$$
\frac{\rho_{n-2}}{\rho_{n-1}}=d_{n},
$$
we also deduce that (\ref{c3t2}) also holds. The definition of $(\tau_n)_n$ also implies that (\ref{c3t3}) holds.
\end{proof}

We finish this section with a couple of comments.

Firstly, Corollary \ref{cop} is not true if we only assume the monic polynomials $\hat p_n$ to have only real zeros without the orthogonality hypothesis. Indeed, It is enough to consider $\hat p_n(x)=x^n$. Take then any sequence $(\rho_n)_n$ of real numbers with $\rho_n\not=0$, $n\ge 0$. For $K=2$ and $\gamma_1=0$, we have
$$
q_n(x)=x^{n-2}(\rho_nx^2+\gamma_2\rho_{n-2}).
$$
We have then that some of the sets
$$
A_-=\{n:\rho_{n}\rho_{n-2}<0\},\quad A_+=\{n:\rho_{n}\rho_{n-2}>0\}
$$
has to be infinite. By taking either $\gamma_2<0$ if $A_-$ is infinite, or
$\gamma_2>0$ if $A_-$ is finite and $A_+$ is infinite, it is clear that
$q_n$ has two non real zeros for infinitely many $n$'s.

And secondly, even if all the zeros of the polynomial $q_n$, $n\ge n_0$, are real, the zeros of $q_{n+1}$ do not necessarily interlace the zeros of $q_n$. Here it is a counterexample: consider the
polynomials $q_n(x)=\hat L_n(x)+4\hat L_{n-1}(x)$, where $\hat L_n$ denotes the $n$-th monic Laguerre polynomial. As a consequence of Lemma \ref{qor}, we have that $q_n$ has only real zeros for $n\ge 0$. But, since
$q_1(x)=x+3$, $q_2(x)=x^2-2$, the zeros of $q_2$ do not interlace the zero of $q_1$.

However, we have found an interesting necessary and sufficient condition for that interlacing property to happen.
Indeed, using the notation of Lemma \ref{lem2}, given real numbers $B_j$, $0\le j\le K-1$, and $\theta$, we can produce the real numbers $A_j$, $0\le j\le K$, as in (\ref{id1l2}), so that we have for
$$
P_{A,\theta}=\sum_{j=0}^{K}A_jx^{K-j},\quad P_{B}(x)=\sum_{j=0}^{K-1}B_jx^{K-1-j}
$$
the identity $P_{A,\theta}(x)=(x-\theta)P_{B}$
(we have included the real number $\theta$ in the notation $P_{A,\theta}$ to stress the dependence of this polynomial on $\theta$).

If we define
$$
q_n^{B}(x)=\sum_{j=0}^{K-1}B_jp_{n-j},\quad q_n^{A;\theta} (x)=\sum_{j=0}^{K}A_jp_{n-j}
$$
the identity (\ref{id1l2}) straightforwardly gives
\begin{equation}\label{idco}
q_{n+1}^{A;\theta}(x)=q_{n+1}^{B}(x)-\theta q_{n}^{B}(x), \quad n\ge K.
\end{equation}
We then have.

\begin{Cor}\label{enze}
Assume that all the zeros of the polynomials $q_n^{B}$ are real and simple for $n\ge n_0$. Then the following conditions are equivalent.
\begin{enumerate}
\item The zeros of $q_{n+1}^{B}$ interlace the zeros of $q_n^{B}$ for $n\ge n_0$.
\item For all real number $\theta$ the polynomial $q_n^{A;\theta}$ has only real and simple zeros for $n\ge n_0+1$.
\end{enumerate}
Moreover, in that case the zeros of $q_n^A$ interlace the zeros of $q_n^B$.
\end{Cor}

\begin{proof}
Using (\ref{idco}), the corollary is an easy consequence of Obreshkov Theorem (see Theorem \ref{obre} in Preliminaries). The interlacing property of the zeros of $q_n^A$ and $q_n^B$ straightforwardly follows from (\ref{idco}).
\end{proof}

As a consequence we have.

\begin{Cor}Let $(p_n)_n$ be a sequence of orthonormal polynomials with respect to a positive measure. Then, for real numbers $\theta_1$, $\theta_2$, the polynomials
$$
q_n(x)=p_n(x)-(\theta_1+\theta_2)p_{n-1}(x)+\theta_1\theta_2p_{n-2}(x)
$$
have only real zeros for $n\ge 0$.
\end{Cor}

\begin{proof}
Since any combination $\tilde q_n(x)=p_n(x)-\theta_1 p_{n-1}(x)$ of orthonormal polynomials has $n$ real zeros, and the zeros of $\tilde q_{n+1}$ interlace the zeros of $\tilde q_n$ (see \cite[Th. 1.2.2]{akh}), it is enough to apply the Corollary \ref{enze}.
\end{proof}

\section{Examples}

We will particularize Theorem \ref{desco} for the classical families of Hermite, Laguerre and Jacobi polynomials.

\subsection{The Hermite case}\label{secH}
In the case of the monic Hermite polynomials $(\hat H_n)_n$ we have that
$$
b_n=0,\quad \hat c_n=n/2
$$
(see \cite[pp. 250-253]{KLS}). As for the zeros we have
$$
\varsigma_{s}=\max\{|\zeta_{j,s}|:1\le j\le s\}\le \sqrt{2s+3}
$$
(see \cite[p.130]{Sze} o also \cite{Ism}).

Let $(\rho_n)_n$ be a normalization sequence of positive real numbers such that
\begin{equation}\label{lamc}
\mbox{$\displaystyle\frac{2\rho_{n-1}}{\rho_n}\le n^\nu$, $n\ge 1$,\quad  for \quad  $0\le \nu\le 1/4$.}
\end{equation}
Write then
$$
\tau_n=\frac{c}{2\sqrt{n+1}},
$$
where $c>0$ will be fixed later on.
We have for $n\ge 2$
\begin{align*}
\max\left\{\left|\frac{\rho_{n-2}}{\rho_{n}\hat c_{n-1}}\right|,\left|\frac{\rho_{n-2}b_{n-1}}{\rho_{n-1}\hat c_{n-1}}\right|,\left|\frac{\rho_{n-2}}{\rho_{n-1}\hat c_{n-1}}\right|\right\}&\le
\max\left\{\frac{(n(n-1))^\nu}{2(n-1)},\frac{(n-1)^\nu}{n-1}\right\},\\
\max\{|\zeta_{j,s}|:1\le j\le s, n-1\le s\le 2n\}&\le \frac{c}{\tau_n}.
\end{align*}
Since the following sequences
$$
\frac{(n(n-1))^\nu\sqrt{n+1}}{2(n-1)},\quad \frac{(n-1)^\nu\sqrt{n+1}}{n-1}
$$
are decreasing in $n$
a simple computation gives that
$$
\max\left\{\left|\frac{\rho_{n-2}}{\rho_{n}\hat c_{n-1}}\right|,\left|\frac{\rho_{n-2}b_{n-1}}{\rho_{n-1}\hat c_{n-1}}\right|,\left|\frac{\rho_{n-2}}{\rho_{n-1}\hat c_{n-1}}\right|\right\}\le \frac{6^\nu}{c}\tau_n,\quad n\ge 3.
$$
If we take $c=6^\nu$ and apply Theorem \ref{desco} to the monic Hermite polynomials, we get the following result.

\begin{Cor}\label{1her} Assume that the normalization sequence $(\rho_n)_n$ satisfies (\ref{lamc}).
For any positive integer $K$ and any finite set of $K+1$ real numbers $\gamma_j$, $j=0,\dots, K$, with $\gamma_0=1$ and $\gamma_K\not=0$, the polynomial
\begin{equation}\label{polhn}
q_n(x)=\sum_{j=0}^K\gamma_j\rho_{n-j}\hat H_{n-j}(x)
\end{equation}
has only real zeros for $n\ge\max\left\{6^{2\nu}\left(\frac{6^{\nu (K-1)}-1}{6^\nu-1}\right)^24^{K-2}\max ^2\{|\gamma_j|,2\le j\le K\},2K\right\}$. Moreover, the zeros are simple and interlace the zeros of $H_{n-1}$.
\end{Cor}

\begin{Rem}
The Hermite family is the case $\rho_n=2^n$, and hence, we can take $\nu=0$. As a consequence, we deduce that all the zeros of the polynomials
$$
q_n(x)=\sum_{j=0}^K\gamma_jH_{n-j}(x)
$$
are real for $n\ge \max\{(K-1)^24^{K-2}\max^2\{|\gamma_j|,2\le j\le K\},2K\}$.
\end{Rem}
Notice that when $(K-1)^24^{K-2}\max^2\{|\gamma_j|,2\le j\le K\}\le 2K$, we have that $q_n$ has only real zeros for $n\ge 2K$.

\bigskip

The probabilistic Hermite polynomials are defined by
$$
\He_n(x)=2^{-n/2}H_n\left(\frac{x}{\sqrt 2}\right).
$$
Hence proceeding as before we conclude that all the zeros of the polynomials
$$
q_n(x)=\sum_{j=0}^K\gamma_j\He_{n-j}(x)
$$
are real for $n\ge\max\{\frac12(K-1)^24^{K-2}\max^2\{|\gamma_j|,2\le j\le K\},2K,11\}$.

\subsection{The Laguerre case}\label{secl}
For $\alpha>-1$, consider the monic Laguerre polynomials $\hat L_n^\alpha$, $n\ge 0$. We
have that
$$
b_n=(2n+\alpha+1),\quad \hat c_n=n(n+\alpha)
$$
(see (\cite[pp, 241-244]{KLS})). As for the zeros we have
$$
\varsigma_{s}=\max\{|\zeta_{j,s}|:1\le j\le s\}\le 4s+2|\alpha|+3
$$
(see \cite[Theorem 6.31.2]{Sze}).

Let $(\rho_n)_n$ be a normalization sequence of positive real numbers such that
\begin{equation}\label{lamcs}
\mbox{$\displaystyle\frac{\rho_{n-1}}{\rho_{n}}\le \nu$, $n\ge 1$,\quad  for \quad  $0< \nu$.}
\end{equation}
Setting $\tau_n=c/(8n+2|\alpha|+3)$, where $c>0$ will be fixed later on, we have for $n\ge (\nu-\alpha+1)/2$
\begin{align*}
\max\left\{\left|\frac{\rho_{n-2}}{\rho_{n}\hat c_{n-1}}\right|,\left|\frac{\rho_{n-2}b_{n-1}}{\rho_{n-1}\hat c_{n-1}}\right|,\left|\frac{\rho_{n-2}}{\rho_{n-1}\hat c_{n-1}}\right|\right\}&\le
\max\left\{\frac{\nu (2n+\alpha-1)}{(n-1)(n-1+\alpha)}\right\},\\
\max\{|\zeta_{j,s}|:1\le j\le s, n-1\le s\le 2n\}&\le \frac{c}{\tau_n}.
\end{align*}
Since the function
$$
f(n,\alpha)=\frac{(2n+\alpha-1)(8n+2|\alpha|+3)}{(n-1)(n-1+\alpha)}, \quad n\ge 2,
$$
is decreasing in $n$ for $\alpha>-1$, we have (since $\alpha+9>2$)
$$
f(n,\alpha)\le f(\alpha+9,\alpha),\quad n\ge \alpha+9.
$$
In turns, a simple computation shows that $f(\alpha+9,\alpha)$ is decreasing for $\alpha>-1$ and
$$
0<f(\alpha+9,\alpha)<23,\quad \alpha>-1.
$$
Hence, taking $c=23\nu$, we have
$$
\max\left\{\left|\frac{\rho_{n-2}}{\rho_{n}\hat c_{n-1}}\right|,\left|\frac{\rho_{n-2}b_{n-1}}{\rho_{n-1}\hat c_{n-1}}\right|,\left|\frac{\rho_{n-2}}{\rho_{n-1}\hat c_{n-1}}\right|\right\}\le \tau_n,\quad n\ge
\max\{(\nu-\alpha+1)/2,\alpha+9\}.
$$
Applying Theorem \ref{desco}, we get the following result.

\begin{Cor}\label{1lag}
Let $\alpha>-1$ and assume that the normalization sequence $(\rho_n)_n$ satisfies (\ref{lamcs}). Then for any positive integer $K$ and any finite set of $K+1$ real numbers $\gamma_j$, $j=0,\dots, K$, with $\gamma_0=1$ and $\gamma_K\not=0$, the polynomial
\begin{equation}\label{qnml}
q_n(x)=\sum_{j=0}^K \gamma_j\rho_{n-j}\hat L^\alpha_{n-j}(x)
\end{equation}
has only real and simple zeros for
$$
n\ge \max\left\{\frac{1}{8}\left(23\nu 2^{K-1}\frac{(23\nu)^{K-1}-1}{23\nu-1}\Gamma_2-2|\alpha|\right),(\nu-\alpha+1)/2,\alpha+9,2K\right\},
$$
and they interlace the zeros of $L^\alpha_{n-1}$ ($\Gamma_2$ defined by (\ref{dgam})).
Moreover, for $n$ big enough all the zeros are positive.
\end{Cor}

\begin{proof}
We have only to prove that the zeros are positive.
Since $\hat L_n^\alpha(0)=(-1)^n(1+\alpha)_n$, we have
$$
q_n(0)=\sum_{j=0}^K(-1)^{n-j}(1+\alpha)_{n-j}\gamma_j.
$$
Hence, we deduce that for $n$ big enough
$$
\sign (q_n(0))=(-1)^n.
$$
Since the zeros of the Laguerre polynomials are positive and the zeros of $q_n$ interlace the zeros of $L^\alpha_{n-1}$, we conclude that for $n$ big enough all the zeros of $q_n$ are positive.
\end{proof}

\begin{Rem}
The monic Laguerre family is the case $\rho_n=1$ ($\nu=1$) and hence, we conclude that all the zeros of the polynomials
$$
q_n(x)=\sum_{j=0}^K\gamma_j\hat L^\alpha _{n-j}(x)
$$
are real for $n\ge \max\{\frac{1}{8}( \frac{23}{22}46^{K-1}\max\{|\gamma_j|,2\le j\le K\}-2|\alpha|),\alpha+9,2K\}$.
\end{Rem}

In the subsequence papers \cite{Durh} and \cite{Durl}, we prove that the spectral properties of the Hermite and Laguerre polynomials allow to improve Corollaries \ref{1her} and \ref{1lag} for certain normalizations of the Hermite and Laguerre polynomials, which include the Hermite and Laguerre polynomials themselves.
Those results are related to the study of zeros
of the polynomial
$$
q_K(x)=\sum_{j=0}^K\gamma _jp_j(x)
$$
(the particular case $n=K$ in (\ref{dqnv})), where $(p_n)_n$ are certain sequences of orthogonal polynomials.
Iserles, N{\o}rsett and Saff \cite{IsSa,IsNo,INS} were probably the first to point out the key role
of the real zeros of the polynomials $\sum_{j=0}^K \gamma_jx^j$, $\sum_{j=0}^K \gamma_j(x)_j$ and some other variants to prove the real rootedness of the polynomial $q_K$,
obtaining some relevant results for Hermite, Laguerre, Gegenbauer, and other important families of orthogonal polynomials.

\subsection{The Jacobi case}\label{secj}
For $\alpha,\beta>-1$, consider the monic Jacobi polynomials $\hat P_n^{\alpha,\beta}$, $n\ge 0$. We
have that
\begin{align*}
b_n&=\frac{\beta^2-\alpha^2}{(2n+\alpha+\beta)(2n+\alpha+\beta+2)},\\
\hat c_n&=\frac{4n(n+\alpha)(n+\beta)(n+\alpha+\beta)}{(2n+\alpha+\beta-1)(2n+\alpha+\beta)^2(2n+\alpha+\beta+1)}
\end{align*}
(see (\cite[pp, 217]{KLS})). As for the zeros we have
$$
|\varsigma_{s}|=\max\{|\zeta_{j,s}|:1\le j\le s\}\le 1
$$
(see \cite[Theorem 3.3.1]{Sze}).

Let $(\rho_n)_n$ be a normalization sequence of positive real numbers such that
\begin{equation}\label{conj}
\mbox{$\rho_{n-1}/\rho_n\le 1/18$, $n\ge 18$,\quad  and\quad  $\lim_n\rho_{n-1}/\rho_n=0$.}
\end{equation}

A simple computation shows that for $n\ge \alpha+\beta+8$
\begin{align*}
\left|\frac{(2n+\alpha+\beta-3)(2n+\alpha+\beta-2)^2(2n+\alpha+\beta-1)}{4(n-1)(n+\alpha-1)(n+\beta-1)(n+\alpha+\beta-1)}\right|&\le 18,\\
\left|\frac{(\beta^2-\alpha^2)(2n+\alpha+\beta-3)(2n+\alpha+\beta-2)^2(2n+\alpha+\beta-1)}{4(n-1)(n+\alpha-1)(n+\beta-1)(n+\alpha+\beta-1)(2n+\alpha+\beta-2)(2n+\alpha+\beta)}\right|&\le 2.
\end{align*}
Hence, for $\tau_n=18\rho_{n-2}/\rho_{n-1}$, we have for $n\ge \max\{\alpha+\beta+8,18\}$
\begin{align}\label{c3t2hj}
\max\left\{\left|\frac{\rho_{n-2}}{\rho_{n}\hat c_{n-1}}\right|,\left|\frac{\rho_{n-2}b_{n-1}}{\rho_{n-1}\hat c_{n-1}}\right|,\left|\frac{\rho_{n-2}}{\rho_{n-1}\hat c_{n-1}}\right|\right\}&\le
\tau_n,\\\label{c3t3hj}
\max\{|\zeta_{j,s}|:1\le j\le s, n-1\le s\le 2n\}&\le \frac{1}{\tau_n}.
\end{align}

Applying Theorem \ref{desco}, we get the following result (the proof is similar to the ones in the previous subsections and it is omitted).

\begin{Cor}\label{1jac}
Let $\alpha,\beta>-1$, and assume that the normalization sequence of positive numbers $(\rho_n)_n$ satisfies (\ref{conj}).
Then for any positive integer $K$ and any finite set of $K+1$ real numbers $\gamma_j$, $j=0,\dots, K$, with $\gamma_0=1$ and $\gamma_K\not=0$, the polynomial
\begin{equation}\label{qnmlj}
q_n(x)=\sum_{j=0}^K \gamma_j\rho_{n-j} \hat P^{\alpha,\beta}_{n-j}(x)
\end{equation}
has only real and simple zeros for $n\ge n_0$, where $n_0$ is the first positive integer such that for $n\ge n_0$,
$\rho_{n-1}/\rho_{n-2}\ge\max\{18(K-1)2^{K-1}\Gamma_2,\alpha+\beta+8,18,2K\}$, and they interlace the zeros of $P^{\alpha,\beta}_{n-1}$. Moreover, for $n$ big enough all the zeros live in $(-1,1)$.
\end{Cor}

\medskip

By taking $\rho_n=n!$, we conclude that all the zeros of the polynomials
$$
q_n(x)=\sum_{j=0}^K\gamma_j n!P_n^{\alpha,\beta}(x)
$$
are real for $n-1\ge\max\{18(K-1)2^{K-1}\Gamma_2,\alpha+\beta+8,18,2K\}$.

\subsection{Asymptotics for the zeros}

The size of the zeros of $q_n$ (\ref{dqnv}) can be deduced taking into account that they interlace the zeros of $p_{n-1}$. For the Hermite, Laguerre and Jacobi cases we also have some asymptotic which we include in the following Corollaries.

\begin{Cor}\label{cas1} Under the hypothesis of Corollary \ref{1her}, for $n$ big enough write $\xi_j(n)$, $1\le j\le n$, for the real zeros of the polynomial $q_n$ (\ref{polhn}) arranged the zeros
in increasing order.
If $n$ is even, then $q_n$ has $n/2$ positive and $n/2$ negative zeros. If $n$ is odd, then $q_n$ has at least $(n-1)/2$ positive and $(n-1)/2$ negative zeros. Moreover, for $k\in\ZZ$,
\begin{align}\label{lec3}
\lim_n 2\sqrt n\xi_{k+n+1}(2n)&=\frac{\pi}2+k\pi,\\\label{lec4}
\lim_n 2\sqrt n\xi_{k+n+1}(2n+1)&=k\pi.
\end{align}
For a bounded continuous function $f$ in $\RR$, we also have
\begin{equation}\label{part2}
\lim_n\frac{1}{n}\sum_{j=1}^nf\left(\frac{\xi_j(n)}{\sqrt{2n}}\right)=\frac{2}{\pi}\int_{-1}^1 f(x)\sqrt{1-x^2}dx.
\end{equation}
\end{Cor}

\begin{proof}
According to Corollary \ref{1her} the zeros of $q_n$ interlace the zeros of $H_{n-1}$. Hence, if $n$ is even $q_n$ has $n/2$ positive and $n/2$ negative zeros, and
if $n$ is odd, then $q_n$ has at least $(n-1)/2$ positive and $(n-1)/2$ negative zeros.

Consider next the well-known Mehler-Heine formula for the Hermite polynomials:
\begin{align}\label{lec}
\lim_n\frac{(-1)^n\sqrt{n\pi}}{2^{2n}n!}H_{2n}\left(\frac{x}{2\sqrt n}\right)&=\cos x,\\\label{lec0}
\lim_n\frac{(-1)^n\sqrt{\pi}}{2^{2n+1}n!}H_{2n+1}\left(\frac{x}{2\sqrt n}\right)&=\sin x
\end{align}
(see \cite[Identities 18.11.7 and 18.11.8]{nist}). This gives after easy computations
\begin{align}\label{lec1}
\lim_n\frac{(-1)^n2^{2n}\sqrt{n\pi}}{2^{2n}n!\rho_{2n}}q_{2n}\left(\frac{x}{2\sqrt n}\right)&=\cos x,\\\label{lec2}
\lim_n\frac{(-1)^n2^{2n+1}\sqrt{\pi}}{2^{2n+1}n!\rho_{2n+1}}q_{2n+1}\left(\frac{x}{2\sqrt n}\right)&=\sin x
\end{align}

The proof of the asymptotic for the central zeros of $q_n$ (\ref{lec3}) and (\ref{lec4}) now follows easily from Hurwitz's Theorem (complex analysis).

In order to prove the limit (\ref{part2}), it is enough to take into account again that the zeros of $q_n$ interlace the zeros of $H_{n-1}$ and
the well-known weak scaling limit
for the counting measure of the zeros of the Hermite polynomials (\cite{Dei})
\begin{equation}\label{wsl}
\lim_n\frac{1}{n}\sum_{j=1}^nf\left(\frac{\zeta_j(n)}{\sqrt{2n}}\right)=\frac{2}{\pi}\int_{-1}^1 f(x)\sqrt{1-x^2}dx,
\end{equation}
where  $\zeta_j(n)$ denote the $j$-th zero of the Hermite polynomial $H_n$.

For the sake of completeness, we describe the exact number of positive and negative zeros of $q_n$ when $n$ is odd, for the Hermite case $\rho_n=2^n$. In order to do that, denote by $\chi$ the following real number: if
$$
\max \{|\gamma_{2j+1}|: 0\le j\le [K/2]\}\not =0,
$$
write $\chi=(-1)^{i+1}\gamma_{2i+1}$, where $2i+1$ is the smallest odd positive integer such that $\gamma_{2i+1}\not =0$, otherwise $\chi=0$.
Then if $n$ is odd, $q_n$ has at least $(n-1)/2$ positive and $(n-1)/2$ negative zeros, plus another zero which is positive if $\chi>0$, negative if $\chi<0$ and zero if $\chi=0$.
Indeed, if $\chi=0$, then for all $1\le j\le [K/2]$, $\gamma_{2j+1}=0$ and then for $n$ odd, $q_n$ is an odd polynomial, hence $q_n(0)=0$. Assume next that $\chi\not =0$.
Since the zeros of $q_n$ interlace the zeros of $H_{n-1}$, we conclude that $q_n$ has at least $(n-1)/2$ positive and $(n-1)/2$ negative zeros, plus another zero in the interval $(\zeta_{\frac{n-1}2}(n-1),\zeta_{\frac{n+1}2}(n-1))$, and
$\zeta_{\frac{n-1}2}(n-1)<0<\zeta_{\frac{n+1}2}(n-1)$. On the one hand, according to the proof of Theorem \ref{desco}, we have
$$
\sign q_n( \zeta_{\frac{n-1}2}(n-1))=\sign H_{n}(\zeta_{\frac{n-1}2}(n-1))=(-1)^{(n-1)/2+1}.
$$
And on the other hand
\begin{align*}
q_n(0)&=\sum_{j=0}^K\gamma_jH_{n-j}(0)=\sum_{j=i}^{[K/2]}(-1)^{(n-2j-1)/2}\frac{(n-2j-1)!}{((n-2j-1)/2)!}\gamma_{2j+1}\\
&=(-1)^{(n-1)/2}\frac{(n-2i-1)!}{((n-1)/2-i)!}\sum_{j=i}^{[K/2]}(-1)^{j}\frac{\left(\frac{n-1}2-i\right)\cdots \left(\frac{n-1}2-j+1\right)}{(n-2i-1)\cdots (n-2j)}\gamma_{2j+1}.
\end{align*}
Hence, for $n$ big enough we have
$$
\sign q_n(0)=(-1)^{(n-1)/2+i}\sign \gamma_{2i+1}=(-1)^{(n-1)/2+1}\chi.
$$
Hence if $\chi<0$, $q_n$ as one more negative zero (in the interval $(\zeta_{\frac{n-1}2}(n-1),0)$). And if $\chi>0$, $q_n$ as one more positive zero (in the interval $(0,\zeta_{\frac{n+1}2}(n-1))$).
\end{proof}

For the Laguerre case we can proceed similarly, using the well-known Mehler-Heine type formula for the Laguerre polynomials (see \cite[Theorem 8.1.3]{Sze}):
\begin{equation*}\label{lel}
\lim_n\frac{L_n^\alpha(z/(n+j))}{n^\alpha}=z^{\alpha/2}J_\alpha(2\sqrt z)=\frac{1}{\Gamma (\alpha+1)}\pFq{0}{1}{-}{\alpha+1}{-z}
\end{equation*}
uniformly in compact set of the complex plane, and the
weak scaling limit for the counting measure of the zeros of the Laguerre polynomials (\cite[Theorem 1]{Gaw})
\begin{equation*}\label{wsl2}
\lim_n\frac{1}{n}\sum_{j=1}^nf\left(\frac{\zeta_j^\alpha(n)}{n}\right)=\frac{1}{2\pi}\int_{0}^4 f(x)\sqrt{\frac{4-x}x}dx,
\end{equation*}
where  $\zeta_j^\alpha(n)$ denote the $j$-th zero of the Laguerre polynomial $L^\alpha_n$.

\begin{Cor}\label{cas1l} Under the hypothesis of Corollary \ref{1lag}, for $n$ big enough write $\xi_j(n)$, $1\le j\le n$, for the real zeros of the polynomial $q_n$ (\ref{qnml}) arranged the zeros
in increasing order. Then, for $i\ge 1$,
\begin{equation*}\label{lel2}
\lim_n n\xi_{i}(n)=j_{i,\alpha}^2/4,
\end{equation*}
where $j_{i,\alpha}$ denotes the $i$-th positive zero of the Bessel function $J_\alpha$.

For a bounded continuous function $f$ in $\RR$, we also have
$$
\lim_n\frac{1}{n}\sum_{j=1}^nf\left(\frac{\xi_j(n)}{n}\right)=\frac{1}{2\pi}\int_{0}^4 f(x)\sqrt{\frac{4-x}x}dx.
$$
\end{Cor}

Finally, for the Jacobi case we can proceed similarly, using the well-known Mehler-Heine type formula for the Jacobi polynomials:
\begin{equation*}\label{lelj}
\lim_n\frac{P_n^{\alpha,\beta}(\cos (z/n))}{n^\alpha}=(z/2)^{-\alpha}J_\alpha(z)
\end{equation*}
uniformly in compact set of the complex plane (see \cite[p. 192]{Sze}), and the
weak scaling limit for the counting measure of the zeros of the Jacobi polynomials (\cite[Section 4.2]{KuVa})
\begin{equation*}\label{wsl2j}
\lim_n\frac{1}{n}\sum_{j=1}^nf\left(\zeta_j^{\alpha,\beta}(n)\right)=\frac{1}{\pi}\int_{-1}^1 \frac{f(x)}{\sqrt{1-x^2}}dx,
\end{equation*}
where  $\zeta_j^{\alpha,\beta}(n)$ denote the $j$-th zero of the Jacobi polynomial $P^{\alpha,\beta}_n$.

\begin{Cor}\label{cas1lj} Under the hypothesis of Corollary \ref{1jac}, for $n$ big enough write $\xi_j(n)$, $1\le j\le n$, for the real zeros of the polynomial $q_n$ (\ref{qnmlj}) arranged the zeros
in increasing order. Then, for $i\ge 1$,
\begin{equation*}\label{lel2j}
\lim_n n\arccos(\xi_{n+1-i}(n))=j_{i,\alpha}.
\end{equation*}
For a bounded continuous function $f$ in $\RR$, we also have
$$
\lim_n\frac{1}{n}\sum_{j=1}^nf\left(\xi_j(n)\right)=\frac{1}{\pi}\int_{-1}^1 \frac{f(x)}{\sqrt{1-x^2}}dx.
$$
\end{Cor}



\end{document}